\documentclass[twoside,fleqn,letterpaper,11pt]{elsarticle}

\usepackage[letterpaper,margin=1in,headheight=0.5in, includefoot]{geometry}
\usepackage{natbib}
\setcitestyle{square,comma}
\usepackage{subcaption}
\usepackage{nomencl} 
\usepackage{booktabs}
\makenomenclature

\usepackage{calligra,amsmath,amssymb}

\def\R{{\mathbb R}}

\newcommand\bfu{{\mathbf u}}

\newcommand\bfI{{\mathbf I}}

\newcommand\bfF{{\mathbf F}}

\newcommand\bfK{{\mathbf K}}
\newcommand\bfL{{\mathbf L}}

\newcommand\bfP{{\mathbf P}}

\newcommand\bfS{{\mathbf S}}

\newcommand\bfU{{\mathbf U}}
\newcommand\bfV{{\mathbf V}}

\newcommand\bfX{{\mathbf X}}

\newcommand\bfZ{{\mathbf Z}}

\newcommand\bfSigma{{\mathbf \Sigma}}

\def\phi{\varphi}

\usepackage[T1]{fontenc}
\journal{arXiv.org}
\begin{document}
\begin{frontmatter}
\title{A Deterministic Dynamical Low-rank Approach for Charged Particle Transport}

\author[1]{Pia Stammer\corref{cor1}}
\cortext[cor1]{p.k.stammer@tudelft.nl}
\author[1]{Tiberiu Burlacu}
\author[2,3]{Niklas Wahl}
\author[1]{Danny Lathouwers}
\author[4]{Jonas Kusch}

\address[1]{Delft University of Technology, Dept. of Radiation Science and Technology, Delft, Netherlands}
\address[2]{German Cancer Research Center -- DKFZ, Heidelberg, Germany}
\address[3]{Heidelberg Institute for Radiation Oncology (HIRO) and National Center for Radiation Research in Oncology (NCRO), Heidelberg, Germany.}
\address[4]{Norwegian University of Life Sciences, \AA s, Norway}

\begin{abstract} 
Deterministically solving charged particle transport problems at a sufficient spatial and angular resolution is often prohibitively expensive, especially due to their highly forward peaked scattering. We propose a model order reduction approach which evolves the solution on a low-rank manifold in time, making computations feasible at much higher resolutions and reducing the overall run-time and memory footprint.
 
For this, we use a hybrid dynamical low-rank approach based on a collided-uncollided split, i.e., the transport equation is split through a collision source method. Uncollided particles are described using a ray tracer, facilitating the inclusion of boundary conditions and straggling, whereas collided particles are represented using a moment method combined with the dynamical low-rank approximation. Here the energy is treated as a pseudo-time and a rank adaptive integrator is chosen to dynamically adapt the rank in energy. We can reproduce the results of a full-rank reference code at a much lower rank and thus computational cost and memory usage. The solution further achieves comparable accuracy with respect to TOPAS MC as previous deterministic approaches. 
\end{abstract}


\end{frontmatter}

\section{Introduction}\label{sec:1}
Charged particle transport problems are relevant to many nuclear as well as medical applications. For example, dose predictions in proton therapy require the solution of high-dimensional transport equations in heterogeneous media with highly forward-peaked scattering. Despite the complexity of the problem, small dose calculation errors and short computation times are required in practice. This often prohibits the use of exact grid- and moment-based numerical methods.

Previous approaches to tackle this problem include for example local refinement of the angular discretization \cite{kophazi_spaceangle_2015,lathouwers_angular_2019,lathouwers_deterministic_2023}. While this can significantly improve the representation of forward-peaked scattering, the system remains high-dimensional and thus costly to solve. Therefore, in this work we want to explore a model order reduction approach.

The dynamical low-rank approximation (DLRA) is a reduced rank approach for large time-dependent matrix differential equations \cite{koch_dynamical_2007}. It evolves the solution in smaller low-rank factors over time, thus maintaining the solution's low-rank structure while minimizing the
residual. Several robust integrators have been proposed \cite{lubich_projector-splitting_2014,ceruti_unconventional_2022,ceruti_rank-adaptive_2022,ceruti_robust_2024} and the efficiency has been shown in a variety of applications including kinetic theory \cite{kusch_robust_2023,ceruti_rank-adaptive_2022,kusch_stability_2023}. In contrast to previous work in this area we now present a 3D, higher order in space framework capable of resolving the highly forward peaked scattering of protons and representing relevant physical effects such as energy straggling. 

This work serves as a proof-of-concept for realistic deterministic low-rank transport calculations with charged particles.  First, we introduce a model for proton transport based on the linear Boltzmann equation with continuous slowing down. We include scattering, continuous slowing down as well as energy straggling terms for uncollided particles. To employ the DLRA
framework in this setting, we split the equation into an uncollided part, which is solved with a raytracer, and a collided part. We then  formulate the energy dependency of the continuous slowing down equation as a pseudo-time and use dynamical low-rank approximation to update the low-rank factors of the collided solution in energy. We present two test cases of a pencil beam hitting a water box and a heterogeneous medium. Our dynamical low-rank approach is compared against TOPAS MC \cite{perl_topas_2012} and the deterministic approach from \cite{lathouwers_deterministic_2023}.

\section{Transport model for proton therapy} \label{sec:2}

\subsection{The linear Boltzmann equation with continuous slowing down}

We model the particle system using the steady-state continuous slowing down (CSD) approximation \cite{larsen_electron_1997} to the linear Boltzmann equation
\begin{align}\label{eq:rte}
    \boldsymbol\Omega \cdot\nabla \psi + \Sigma_s(E,\mathbf x) \psi = \frac{\partial S\psi}{\partial E} + \frac12 \frac{\partial^2 T\psi}{\partial E^2} + \int_{\mathbb{S}^2} \Sigma_s(E,\mathbf x,\boldsymbol\Omega\cdot\boldsymbol\Omega')\psi(E,\mathbf x,\boldsymbol\Omega')\mathrm{d}\boldsymbol\Omega'.
\end{align}
 The phase space of the particle density $\psi$ consists of energy $E\in [0, E_{max}]\subset\mathbb{R}_+$, space $\mathbf{x}\in \mathcal{X}\subset\mathbb{R}^3$ and direction of flight $\boldsymbol{\Omega}\in\mathbb{S}^2$. $\Sigma_s$ denotes scattering cross sections, which describe the likelihood of (in- or out-)scattering interactions of particles with tissue. Further, we use $S:\mathbb{R}_+ \times\mathbb{R}^3 \rightarrow\mathbb{R}_+$ to denote the stopping power, which describes the expected energy loss and $T$ for the straggling coefficient \cite[see e.g.][]{burlacu_deterministic_2023}.
The tissue density of the patient is $\rho:\mathcal{X}\rightarrow\mathbb{R}_+$. For simplicity, we follow the common assumption that all materials are water-equivalent and differ only in density \cite[e.g.][]{woo_validity_1990, olbrant_generalized_2010}, i.e.,
\begin{align*}
T(E,\mathbf x) = T^{H_2O}(E)\rho(\mathbf x), \quad S(E,\mathbf x) = S^{H_2O}(E)\rho(\mathbf x), \quad \Sigma_s(E,\mathbf x,\mathbf \Omega\cdot\mathbf \Omega') = \Sigma_s^{H_2O}(E,\mathbf \Omega\cdot\mathbf \Omega')\rho(\mathbf x)\;.
\end{align*}

\subsection{Physics data}
The material and particle dependent values for stopping power, straggling and scattering are then determined using physical formulas as well as databases. In the following we neglect nuclear scattering effects since their contribution to the spatial dose distribution is relatively small \cite{newhauser_physics_2015,durante_nuclear_2016}.
\subsubsection{Scattering}
We model the differential Coulomb scattering cross sections in water, $\Sigma_s^{H_2O}(E,\mathbf \Omega\cdot\mathbf \Omega')$, using a screened Moli\`{e}re model \cite[see e.g.][]{burlacu_deterministic_2023,scott_theory_1963} and the corresponding total cross section is obtained by integrating over $\mathbf \Omega \in \mathbb{S}^2$.

\subsubsection{Stopping power}
Values for the stopping power are retrieved using PStar, based on the NIST database \cite{mj_berger_js_coursey_ma_zucker_and_j_chang_stopping-power_2009} or can be computed using Bethe-Bloche's model \cite{bethe_zur_1930,bloch_zur_1933}.
\subsubsection{Straggling}
 While the stopping power models expected energy loss, interactions with the material are stochastic and thus the amount of deposited energy through such interactions differs \cite{bethe_zur_1930,bloch_zur_1933}. 
 We use William's model \cite{williams_passage_1997,burlacu_deterministic_2023} to compute the modified straggling coefficient $T^{H_2O}(E)$ in water for the uncollided part of the equation. 
\section{Methods}
\subsection{Collided-uncollided split}
We wish to solve the evolution equation \eqref{eq:rte}.
 We now separate the angular flux into collided particles $\psi_c$ and uncollided particles $\psi_u$ such that $\psi = \psi_u + \psi_c$. For collided particles, we assume that straggling can be neglected, that is, $ \frac{\partial^2 T\psi_c}{\partial E^2} \approx 0.$
Then, the original equation \eqref{eq:rte} becomes
\begin{align}\label{eq:rte-cu}
    \mathbf{\Omega} \cdot\nabla \psi_u +  \mathbf{\Omega} \cdot\nabla \psi_c + \Sigma_s(E,\mathbf{x}) (\psi_u + \psi_c)=& \frac{\partial S\psi_u}{\partial E} + \frac{\partial S\psi_c}{\partial E} + \frac12 \frac{\partial^2 T\psi_u}{\partial E^2}\nonumber \\ &+ \int_{\mathbb{S}^2} \Sigma_s(E,\mathbf{x}, \mathbf{\Omega}\cdot \mathbf{\Omega'})\left(\psi_u(E,\mathbf{x},\mathbf\Omega') + \psi_c(E,\mathbf x,\mathbf\Omega')\right)\mathrm{d}\mathbf\Omega'
\end{align}
which can be separated into the coupled system
\begin{subequations}\label{eq:rte-cu-split} 
    \begin{align}
        \mathbf\Omega \cdot\nabla \psi_u + \Sigma_s(E,\mathbf x) \psi_u =& \frac{\partial S\psi_u}{\partial E} + \frac12 \frac{\partial^2 T\psi_u}{\partial E^2}, \\ 
        \mathbf\Omega \cdot\nabla \psi_c + \Sigma_s(E,\mathbf x) \psi_c =& \frac{\partial S\psi_c}{\partial E} + \int_{\mathbb{S}^2} \Sigma_s(E,\mathbf x,\mathbf \Omega\cdot \mathbf\Omega')\left(\psi_u(E,\mathbf x,\mathbf\Omega') + \psi_c(E,\mathbf x,\mathbf\Omega')\right)\mathrm{d}\mathbf\Omega'. \label{eq:rte-cu-split-b}
    \end{align}
\end{subequations}

\subsection{Raytracer for the uncollided part}
We use an in-house raytracer with 128 energy groups and a second order polynomial basis in energy to compute the uncollided flux along rays \cite[see also][]{burlacu_deterministic_2023,lathouwers_deterministic_2023}. The (rectangular) spatial domain is discretized using a regular grid and a beam source is placed at the boundary on one side of the domain. 
Rays are traced from the midpoints of the boundary faces of the 3D grid cells in this boundary plane through the domain. In this work we assume an initially unidirectional beam, i.e.,  $\psi_u(E,\mathbf{x},\mathbf{\Omega}') = \psi(E,\mathbf{x})\delta(\mathbf{\Omega}' - \mathbf{\Omega}_{\mathrm{in}})$. Thus, for the uncollided flux, we only have to trace along straight lines. The uncollided flux at each grid cell and energy then contributes to the dynamics of the collided particles through inscattering as described in equation \eqref{eq:rte-cu-split-b}.
\subsection{Dynamical low-rank approximation for the collided part}
We now want to use the dynamical low-rank approximation to solve the collided equation \eqref{eq:rte-cu-split-b} more efficiently. Since the problem is not time dependent, we apply a transformation as described in \cite{kusch_robust_2023} and let the energy act as a pseudo time $t$: 
\begin{align}
    t(E) := \int_{E}^{E_{\mathrm{max}}}\frac{1}{S(E')}\,dE', \quad 
    \widetilde{\psi}(t,\mathbf x,\mathbf  \Omega) = \rho(\mathbf x)S(E(t))\psi(E(t),\mathbf x,\mathbf \Omega), \quad \Sigma_s(t,\mathbf \Omega) = \Sigma_s(E(t),\mathbf \Omega),
\end{align}
which yields
\begin{align}\label{eq:trafo-split-b}
    \partial_t \widetilde{\psi}_c + \mathbf \Omega \cdot\nabla \frac{\widetilde{\psi}_c}{\rho(\mathbf x)} + \Sigma_s(t,\mathbf x) \widetilde{\psi}_c =& \int_{\mathbb{S}^2} \Sigma_s(t,\mathbf x,\Omega\cdot\Omega')\left(\widetilde{\psi}_u(t,\mathbf x,\mathbf \Omega') + \widetilde{\psi}_c(t,\mathbf x,\mathbf \Omega')\right)d\mathbf \Omega'.
\end{align}
\subsubsection{Spatial and angular discretization}
Next, we discretize in space and angle to obtain a (pseudo)-time dependent matrix differential equation.
For the angular discretization we choose the modal  spherical harmonics (P$_N$) method \cite{case_linear_1968} and collect all real-valued spherical harmonics basis functions up to degree $N$ in a vector
\begin{align*}
    \mathbf m = (m_0^0, m_1^{-1}, m_1^{0}, m_1^{1},\cdots, m_N^{N})^{\top}\in\mathbb{R}^{(N+1)^2}.
\end{align*}
This gives the modal approximation $\widetilde{\psi}_c(t,\mathbf{x},\mathbf \Omega) \approx \mathbf u(t,\mathbf x)^{\top}\mathbf m(\mathbf\Omega)$. Then, the P$_N$ equations read
\begin{align*}
    \partial_t \mathbf u (t,\mathbf x) =-\mathbf A\cdot\nabla \frac{\mathbf u(t,\mathbf x)}{\rho(\mathbf x)}-\Sigma_t(t) \mathbf u(t,\mathbf x)+\mathbf G \mathbf u(t, \mathbf x) + \mathbf G \mathbf M\widetilde{\psi}_u (t, \mathbf x),\;
\end{align*}
where $\mathbf A\cdot\nabla := \mathbf A_1\partial_{x} + \mathbf A_2\partial_y+ \mathbf A_3\partial_z$ with $\mathbf A_i := \int_{\mathbb{S}^2} \mathbf m \mathbf m^{\top} \Omega_i \, \mathrm{d}\mathbf\Omega$. The diagonal in-scattering matrix $\mathbf G$ has entries $G_{kk}(t) = 2\pi\int_{[-1,1]}P_{k}(\mu)\Sigma_s(t,\mu)\,d\mu$. We use the matrix $\mathbf M = (w_k m_i(\Omega_k))_{ik}$ which transforms the nodal solution from the ray-tracer into its modal counterpart. 

Finally, we discretize in space using finite volumes with a second order upwind scheme \cite[see also][]{kusch_robust_2023,kusch_kit-rt_2023}.
\subsubsection{Time/energy integration using the dynamical low-rank approximation}
Having discretized in space and angle, our remaining task is to evolve the solution of a large matrix ordinary differential equation
\begin{align}
    \dot \bfu(t) = \bfF(t, \bfu(t))\, , \quad \bfu(t_0) = \bfu_0\,,
\end{align}
where $\bfu\in\R^{n_x\times n_{\Omega}}$. Typically, $n_x$ is the number of spatial cells, and $n_{\Omega}$ is the number of additional parameter points, such as the number of moments of the P$_N$ approximation. Then, a rank $r$ approximation to $\bfu(t)$ can be written as a SVD-like factorization $\bfu(t)\approx \bfX(t)\bfS(t)\bfV(t)^{\top}$, where $\bfX(t)\in\R^{n_x \times r}$ and $\bfV(t)\in\R^{n_{\Omega} \times r}$ are orthonormal basis matrices and $\bfS(t)\in\mathbb{R}^{r\times r}$ is an invertible, but not necessarily diagonal, coefficient matrix. Note that the rank $r$ matrices form a smooth manifold, which we will denote by $\mathcal{M}$. The core idea of dynamical low-rank approximation (DLRA) is to evolve the low-rank factors instead of the full solution in time \cite{koch_dynamical_2007}. This is achieved by projecting the dynamics onto the tangent space of $\mathcal{M}$. The projection ensures that the solution will remain on the low-rank manifold as time evolves. We are then interested in solving
\begin{align}\label{eq:projectedDyn}
    \dot \bfu_r(t) = \bfP(\bfu_r(t))\bfF(t, \bfu_r(t))\,,
\end{align}
where for $\bfZ=\bfX\bfSigma\bfV^{\top}$, the projector onto the tangent space of $\mathcal{M}$ at $\bfZ$ takes the form 
\begin{align*}
    \bfP(\bfZ)\bfF = \bfX\bfX^{\top}\bfF(\bfI-\bfV\bfV^{\top}) + \bfF\bfV\bfV^{\top}\,.
\end{align*}

Since the manifold can exhibit high curvature, conventional time integration methods require prohibitively small step sizes. High curvature arises when the solution exhibits small singular values, typically occurring when the selected rank allows for an accurate approximation of the solution. Consequently, high curvature naturally occurs throughout the entire simulation. Therefore, careful construction of new time integration methods that take the geometry of $\mathcal{M}$ into account when solving \eqref{eq:projectedDyn} is required \cite{kieri_discretized_2016}. Robust time integration methods include the projector--splitting integrator~\cite{lubich_projector-splitting_2014} and basis-update \& Galerkin (BUG) integrators~\cite{ceruti_unconventional_2022,ceruti_rank-adaptive_2022}. The main idea of these integrators is to move on flat subspaces of the low-rank manifold, thereby avoiding regions of high curvature. In this work, we use the rank-adaptive/augmented BUG integrator \cite{ceruti_rank-adaptive_2022}, which first updates the basis matrices $\bf{U}$ and $\bfV$ in parallel, constraining the dynamics to
\begin{align}
    \{\bfK \bfV(t_0)^{\top} | \bfK \in\R^{n_x \times r}\} \subset \mathcal{M} \enskip\text{ and }\enskip \{\bfU(t_0)\bfL^{\top} | \bfL \in\R^{n_{\Omega} \times r}\}\,\subset \mathcal{M}.
\end{align}
This leads to the differential equations
\begin{subequations}\label{eq:KLS}
\begin{align}
    \dot{\bfK}(t) =\,& \bfF(t, \bfK(t)\bfV_0^{\top})\bfV_0\,,\quad \,&&\bfK(t_0) = \bfU_0\bfS_0\,,\\
    \dot{\bfL}(t) =\,& \bfF(t, \bfU_0\bfL(t)^{\top})^{\top}\bfU_0\,,\quad \,&&\bfL(t_0) = \bfV_0\bfS_0^{\top}\,.
\end{align}
The augmented bases $\mathbf{\widehat{U}}\in\R^{n_x \times 2r}$ and $\mathbf{\widehat{V}}\in\R^{n_{\Omega} \times 2r}$ can then be retrieved by a QR decomposition of $[\bfK(t_1) \; \bfU_0]$ and $[\bfL(t_1) \; \bfV_0]$, respectively. To update the augmented coefficient matrix $\mathbf{\widehat{S}} \in\R^{2r \times 2r}$ we solve the Galerkin system
\begin{align}
    \dot{\mathbf{S}}(t) =\,& \mathbf{\widehat{U}}^{\top}\mathbf{F} (t, \mathbf{\widehat{U}}\mathbf{S}( t )\mathbf{\widehat{V}}^{\top})\mathbf{\widehat{V}}\,,\quad \mathbf{S} (t_0) = (\mathbf{\widehat{U}}^{\top}\mathbf{U}_0)\mathbf{S_0}(\mathbf{V}_0^{\top}\mathbf{\widehat{V}})\,.
\end{align}
\end{subequations}
Finally, the time-updated basis and coefficient matrices $\bfU_1,\bfV_1,\bfS_1$ at time $t_1 = t_0 + \Delta t$ are obtained by truncating  $\bf\widehat{U}$,  $\bf\widehat{V}$ and  $\bf\widehat{S}$. We follow \cite{ceruti_rank-adaptive_2022} and truncate according to the tolerance $\vartheta$, choosing the new rank $r_{1}\le 2r$ such that
		$$
		\biggl(\ \sum_{j=r_{1}+1}^{2r} \sigma_j^2 \biggr)^{1/2} \le \vartheta\;,
		$$
where $\sigma_j$ are the singular values of  $\bf\widehat{S}$.

The complete process is repeated until the factorized solution at a given final time $t_n$ is reached after $n$ steps.
Note that this algorithm adapts the rank in each time step based on a user-defined tolerance $\vartheta$ which represents the (estimated) truncation error. The rank-adaptivity thereby takes into account and can give insight into the changing dynamics and complexity of the solution at different points in time.\\ 
Solving equations \eqref{eq:KLS} will avoid any curvature from $\mathcal{M}$. Therefore, classical time integration methods can be used. Here, we follow the approach in \cite{kusch_robust_2023} and use an implicit time discretization on the scatter part of the integrator and explicit updates on the remainder. This approach has been shown to be stable under a time-step restriction that is independent of small
singular values in $\mathbf S$ and stiff terms due to scattering \cite{kusch_robust_2023}. 

By construction, the augmented BUG integrator is first order in pseudo-time/energy. Higher rates can however sometimes be observed in practice \cite{ceruti_rank-adaptive_2022,ceruti_robust_2024}. To overcome the limitation to first order, the projector-splitting integrator \cite{lubich_projector-splitting_2014} or recent second-order variations of the BUG integrator \cite{ceruti_robust_2024,kusch_second-order_2024} could be used. The former is however less stable for diffusive problems \cite{kusch_stability_2023}. The latter is associated with higher costs due to an augmentation to $3r$ or $4r$ and has been observed to be less accurate than the augmented BUG integrator \cite{ceruti_robust_2024}.

\section{Numerical results}
To test the accuracy and performance of our method, we consider two test cases in a 2x2x8cm domain irradiated with a unidirectional 90MeV Gaussian proton beam with spatial standard deviation $\sigma = 0.3$cm and a standard deviation of 1MeV in energy. We compare our dynamical low-rank approach at different angular resolutions and tolerances to a TOPAS MC reference with only electro-magnetic interactions and $10^7$ simulated histories. The spatial grid for all methods has a resolution of 0.5mm, i.e., 40x40x160 cells. The deposited energy is normalized to the average per particle in the Monte Carlo reference.
\subsection{Computational costs and memory}
The dynamical low-rank approximation reduces the computational costs of a computation with $n$ spatial cells and $m$ angular degrees of freedom from $O(n\cdot m$) to $O(r^2\cdot(n+m))$ and the required memory to $O(r\cdot (n+m))$ compared to a full-rank P$_N$ approach \cite{kusch_robust_2023}. The runtime on the same machine reduces to less than 1\%. This makes it possible to use angular resolutions that are not previously computationally feasible. 
\subsection{90 MeV homogeneous}
\begin{figure}[h] 
     \centering
     \begin{subfigure}[b]{0.425\textwidth}
         \centering
         \includegraphics[width=\textwidth]{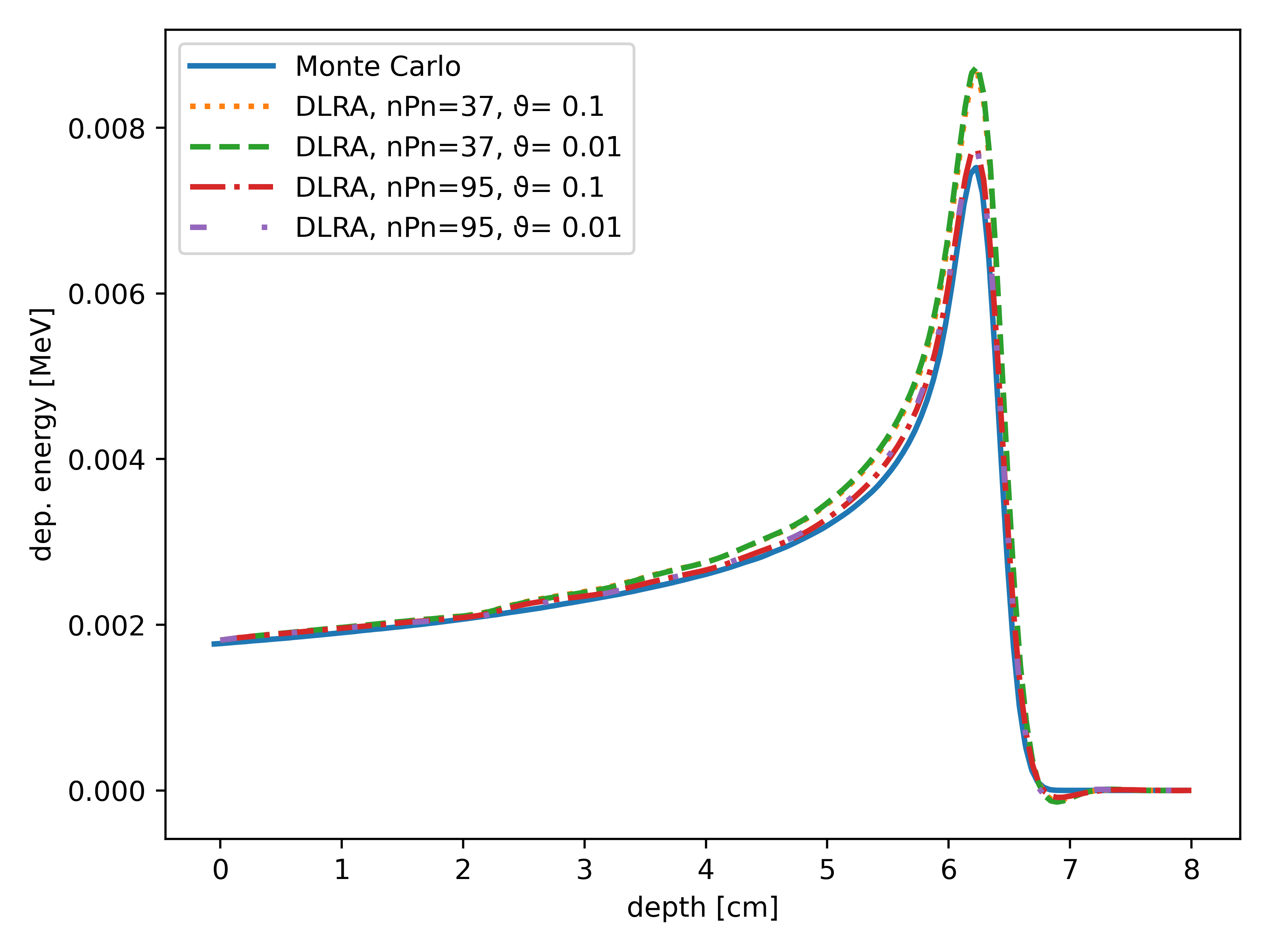}
         \caption{Longitudinal cut at x = y = 1cm.}
     \end{subfigure}
          \begin{subfigure}[b]{0.425\textwidth}
         \centering
         \includegraphics[width=\textwidth]{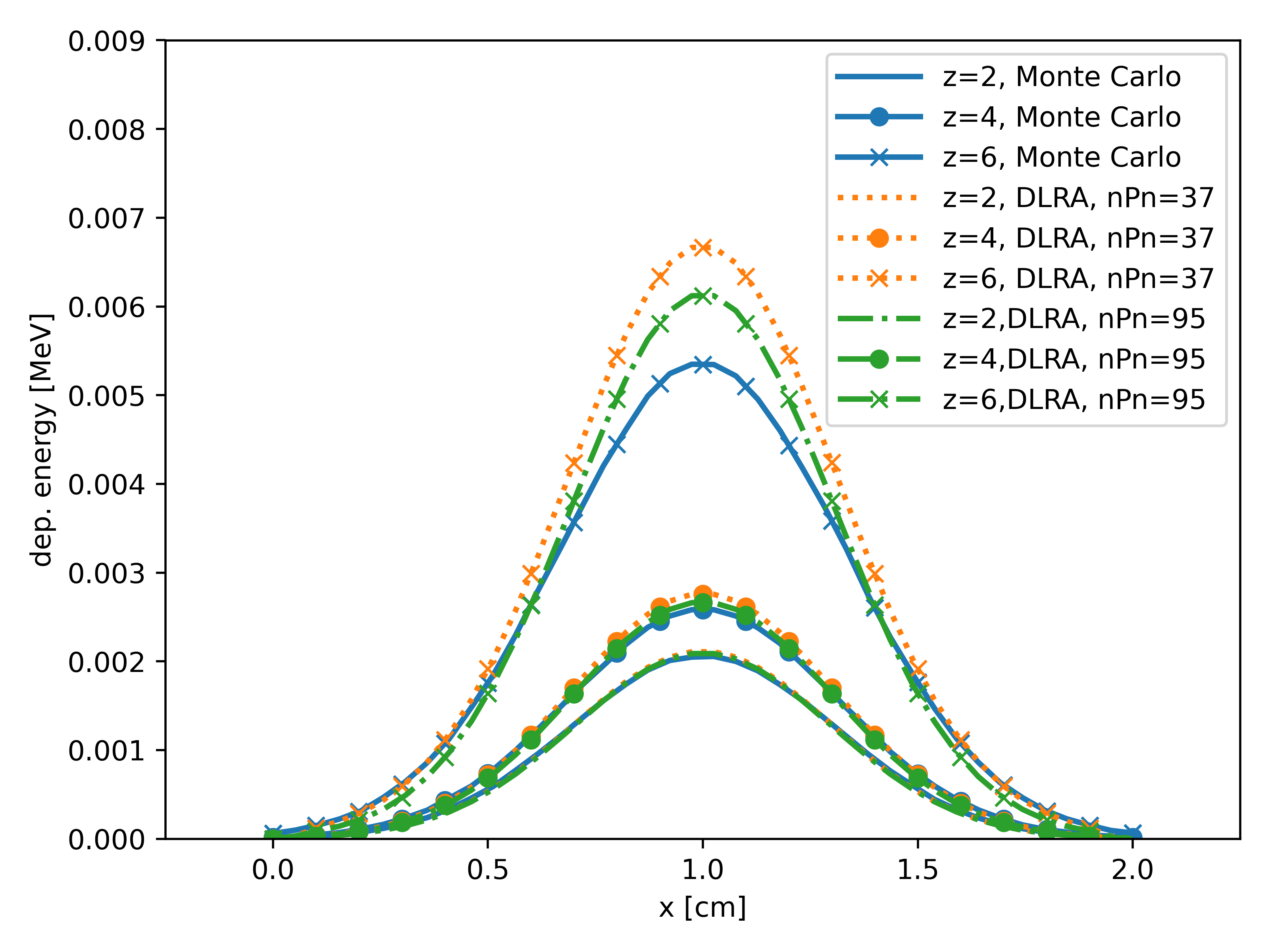}
         \caption{Lateral cuts at y = 1cm for $\vartheta=0.1$.}
     \end{subfigure}
        \caption{Cuts through the deposited energy for Monte Carlo vs. DLRA with different angular resolutions and tolerances for homogeneous case.}
        \label{fig:HomCuts}
\end{figure}

\begin{figure}[h] 
\centering
     \begin{subfigure}[b]{0.2\textwidth}
         \centering
         \includegraphics[width=.75\textwidth]{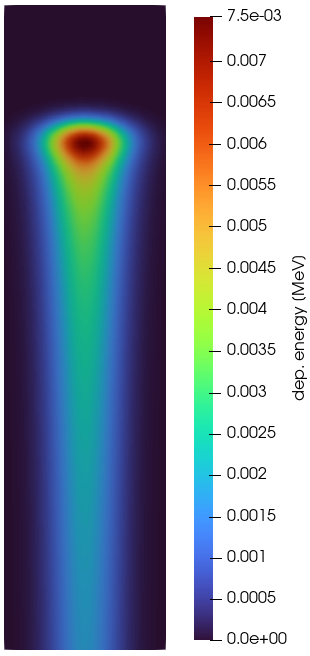}
         \caption{Monte Carlo}
     \end{subfigure}
          \begin{subfigure}[b]{0.2\textwidth}
         \centering
         \includegraphics[width=.75\textwidth]{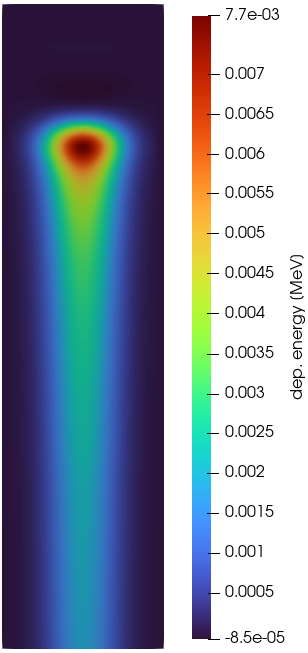}
         \caption{DLRA}
     \end{subfigure}
          \begin{subfigure}[b]{0.2\textwidth}
         \centering
         \includegraphics[width=.75\textwidth]{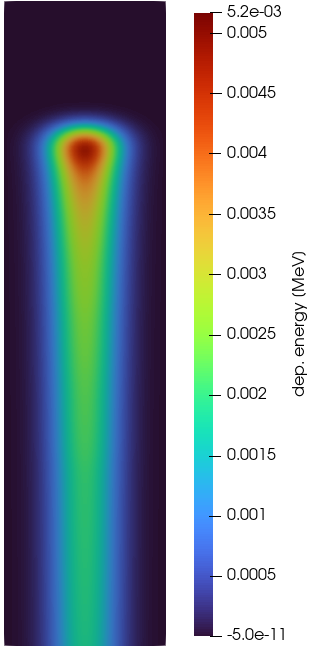}
         \caption{DLRA uncollided}
     \end{subfigure}
          \begin{subfigure}[b]{0.2\textwidth}
         \centering
         \includegraphics[width=.75\textwidth]{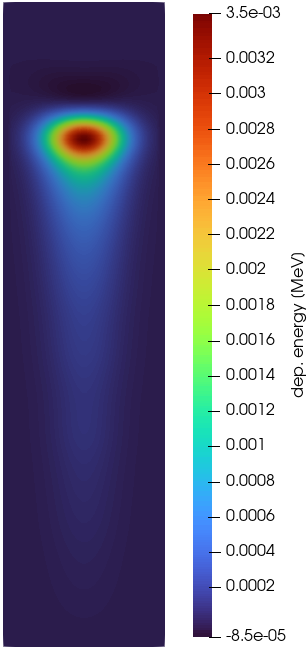}
           \caption{DLRA collided}
     \end{subfigure}
        \caption{Deposited energy of MC vs. DLRA with P$_{95}$, $\vartheta=0.1$ at x=1cm in the homogeneous test case.}
        \label{fig:HomSurfPlots}
\end{figure}

In the first test case, we consider a homogeneous water box with 0 Hounsfield units (HU). Figures \ref{fig:HomCuts} and \ref{fig:HomSurfPlots} show that the low-rank approach can capture the important characteristics of both slowing down and scattering of particles. Overall, the deposited energy agrees well with TOPAS MC for the higher angular resolution apart from small deviations around the Bragg peak and a small amount of negative deposition
 behind the peak. 

\begin{figure}[h] 
     \centering
     \begin{subfigure}[b]{0.3275\textwidth}
         \centering
         \includegraphics[width=\textwidth]{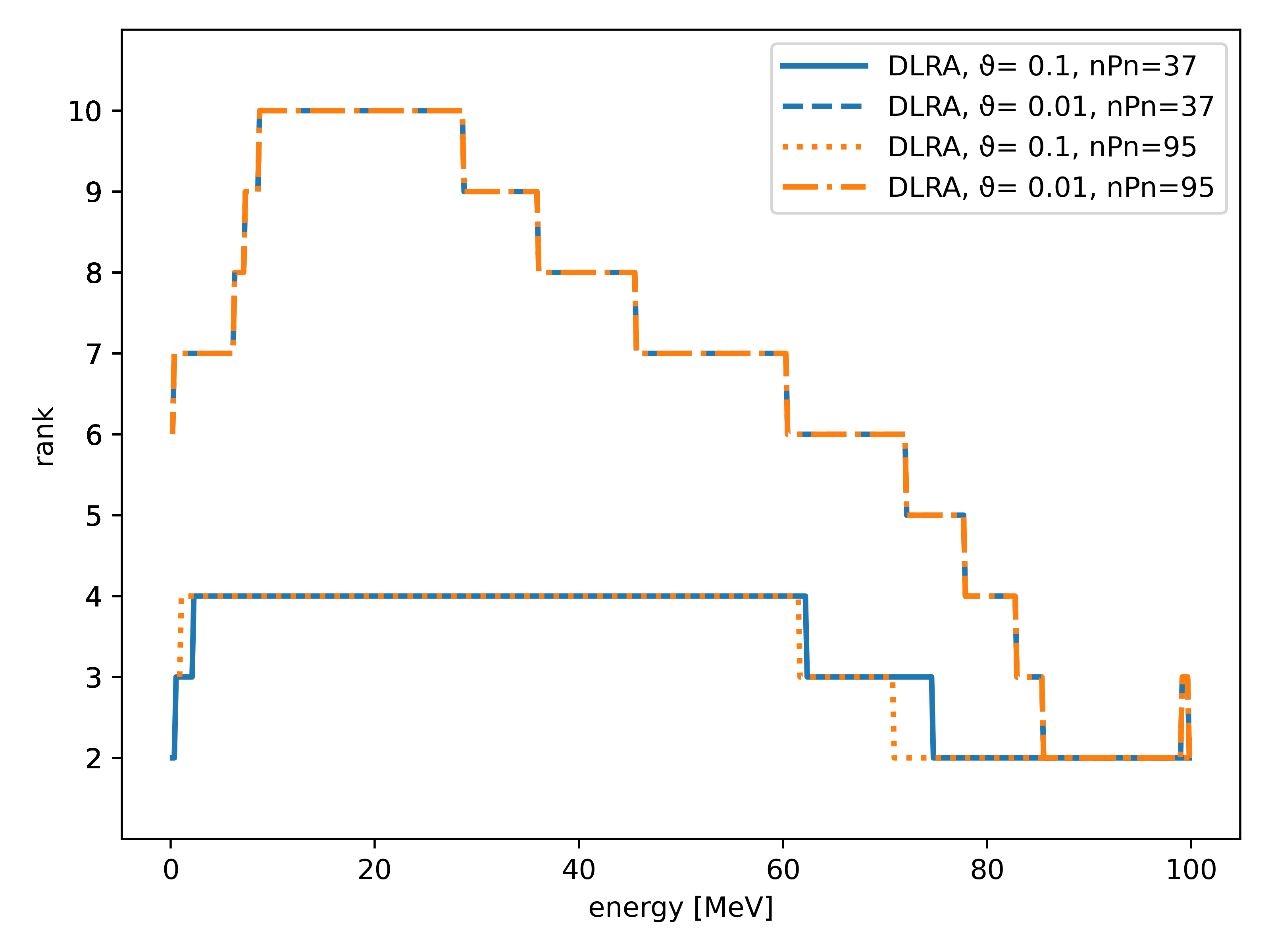}
         \caption{Ranks, homogeneous test case.}
     \end{subfigure}
          \begin{subfigure}[b]{0.3275\textwidth}
         \centering
         \includegraphics[width=\textwidth]{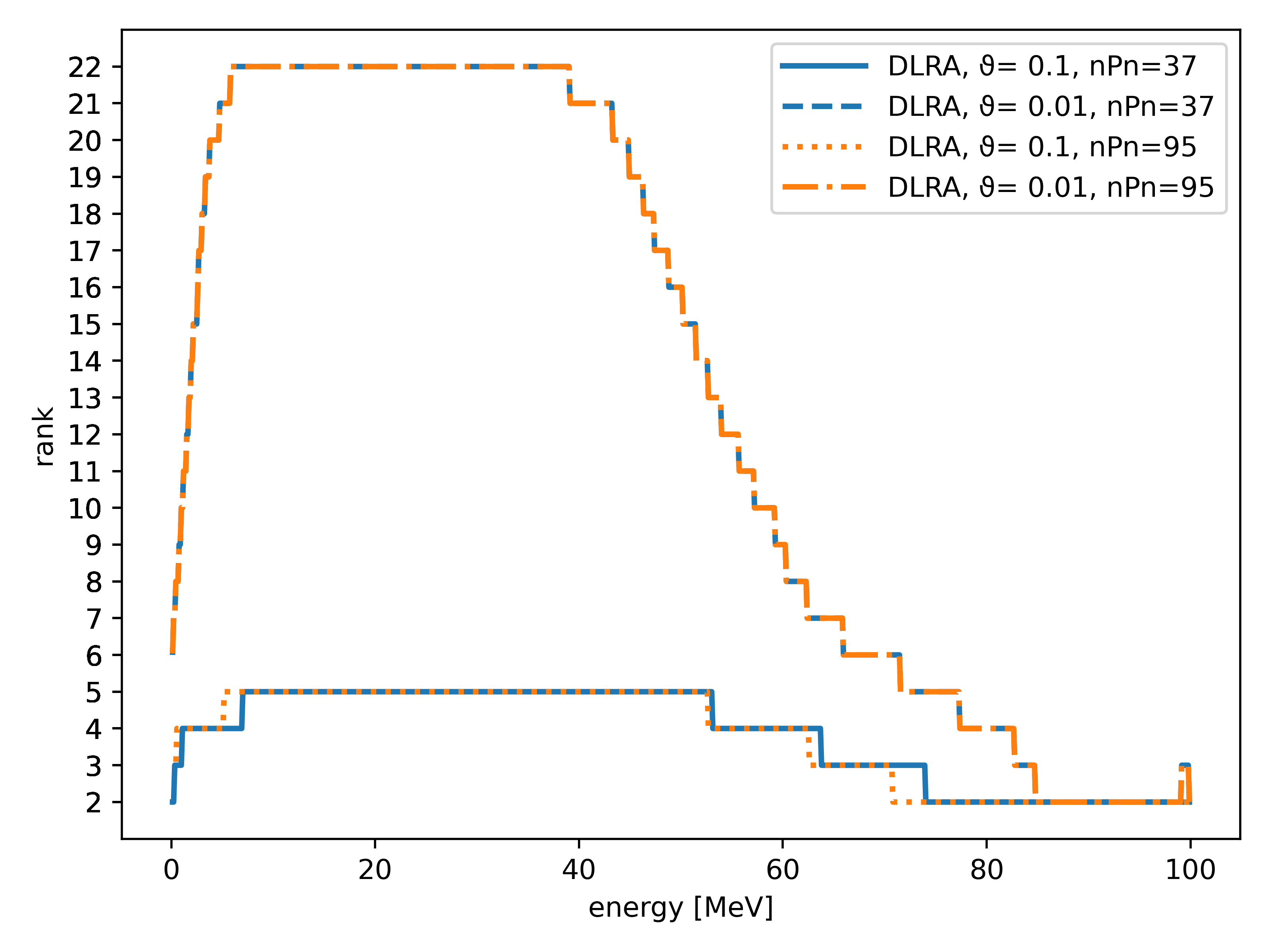}
         \caption{Ranks, heterogeneous test case.}
     \end{subfigure}
        \begin{subfigure}[b]{0.32125\textwidth}
         \centering
         \includegraphics[width=\textwidth]{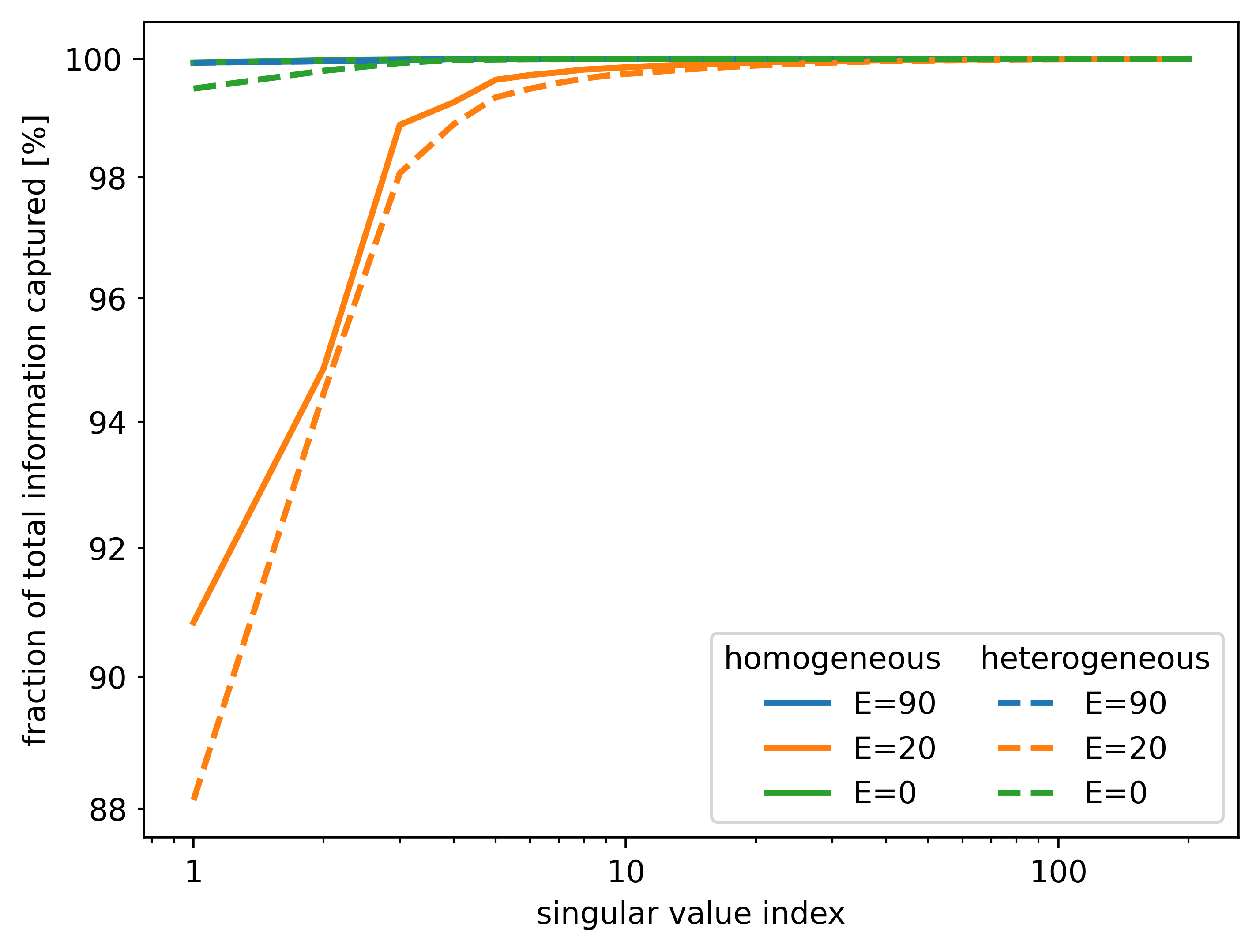}
         \caption{Captured information.}
     \end{subfigure}
        \caption{(a-b) Ranks chosen by rank-adaptive integrator for different tolerances and angular resolutions. (c) Fraction of captured information in top k singular values\protect\footnotemark of P$_{95}$ solution at different energies.}
        \label{fig:Ranks}
\end{figure}

This aligns with the results from \cite{lathouwers_deterministic_2023}, despite the fact that we neglect straggling in the collided part of the equation, indicating that this is a reasonable approximation. Further, negative energy deposition is also observed using the full-rank approach in \cite{lathouwers_deterministic_2023} which implies that the cause are the chosen numerical discretization or model approximations rather than the low-rank approach. While we suffer less from spatial resolution issues compared to \cite{lathouwers_deterministic_2023}, the remaining differences to Monte Carlo could be attributed partly to the limited order of convergence in energy as well as the angular discretization or limitations of our simplified model.
Indeed, figure \ref{fig:HomCuts} (a) as well as table \ref{tab:L2diff} show that a higher number of moments is necessary to sufficiently resolve the collided part in angle. However at a sufficient angular resolution, the choice of the tolerance parameters for the rank adaption does not have a large impact on the solution. Thus, an even higher angular resolution might further improve the DLRA results. Since the required ranks chosen over pseudo time/energy (figure \ref{fig:Ranks}) are low and do not increase significantly from P$_{37}$ to P$_{95}$, this should still be computationally feasible. 
\begin{table}[h]
    \centering
    \caption{Relative discrepancy between Monte Carlo and DLRA in $L^2$ norm.}
    \begin{tabular}{ccccc}
    \toprule 
        &P$_{37}$, $\vartheta =0.1$ &P$_{37}$, $\vartheta =0.01$ &P$_{95}$, $
        \vartheta =0.1$ & P$_{95}$, $\vartheta =0.01$\\
        \midrule
        homogeneous &  0.058 & 0.047 & 0.041	&0.039 \\
        heterogeneous & 0.066 & 0.067 &	0.059	& 0.057 \\
        \bottomrule
    \end{tabular}
    \label{tab:L2diff}
\end{table}
However, the fact that the ranks do not increase for finer angular resolutions also indicates that a smaller number of angular modes can be sufficient if chosen appropriately. Thus, a combination of the low-rank approach with local angular refinement could be beneficial.

\subsection{90 MeV heterogeneous}
\begin{figure}[h!] 
     \centering
     \begin{subfigure}[b]{0.425\textwidth}
         \centering
         \includegraphics[width=\textwidth]{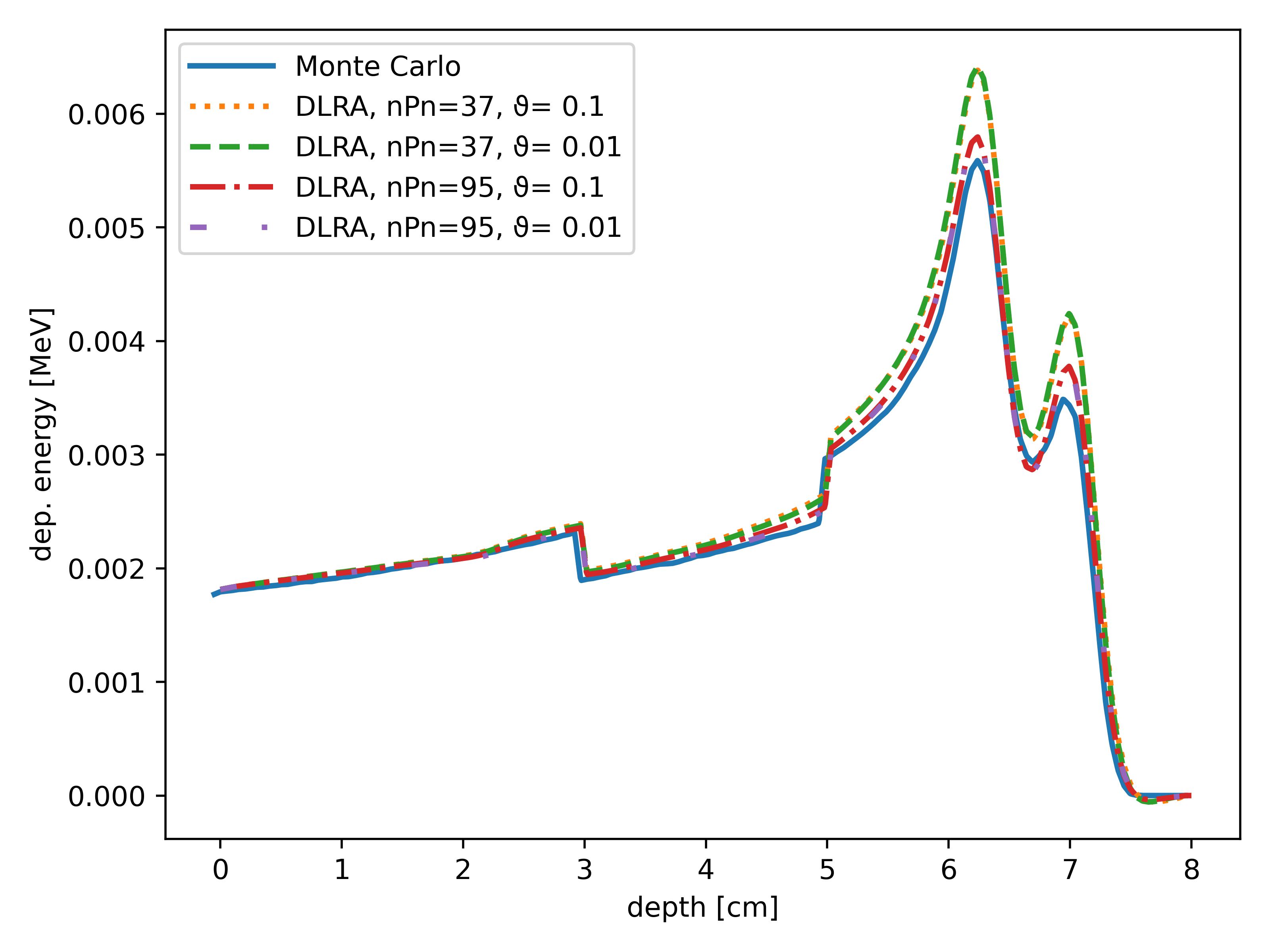}
         \caption{Longitudinal cut at x = y = 1cm.}
     \end{subfigure}
          \begin{subfigure}[b]{0.425\textwidth}
         \centering
         \includegraphics[width=\textwidth]{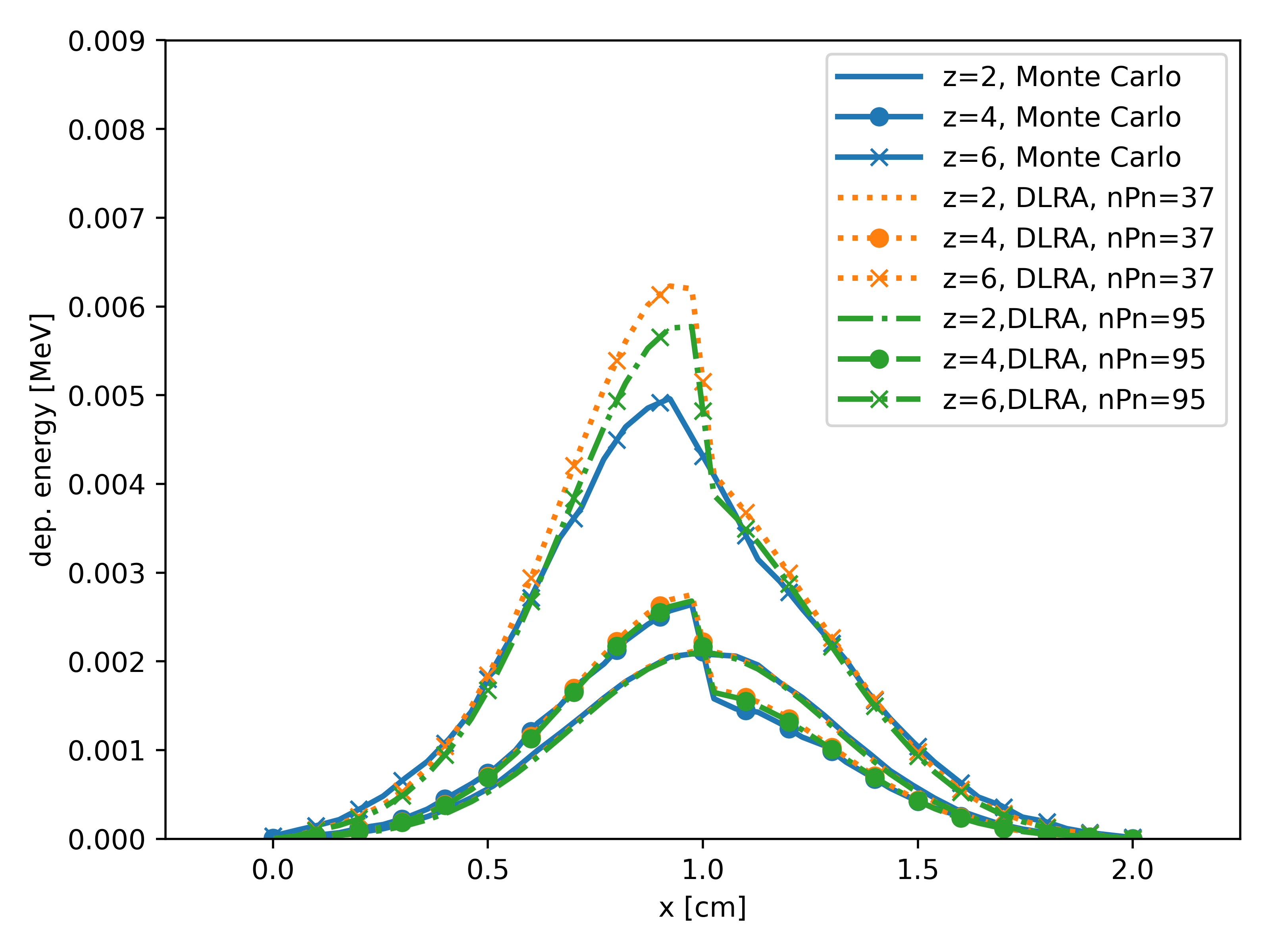}
         \caption{Lateral cuts at y = 1cm for $\vartheta=0.1$.}
     \end{subfigure}
        \caption{Cuts through the deposited energy for Monte Carlo vs. DLRA with different angular resolutions and tolerances in heterogeneous case.}
        \label{fig:HetCuts}
\end{figure}

\begin{figure}[h!] 
     \centering
     \begin{subfigure}[b]{0.2\textwidth}
         \centering
         \includegraphics[width=.75\textwidth]{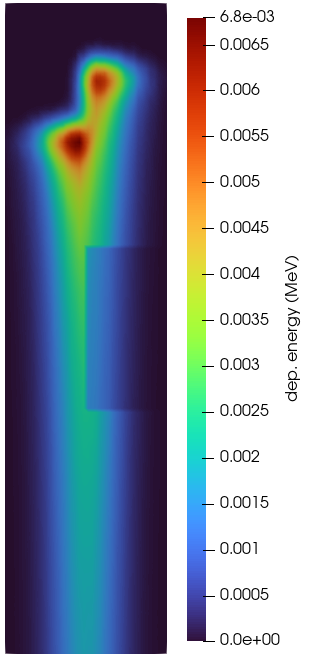}
     \end{subfigure}
          \begin{subfigure}[b]{0.2\textwidth}
         \centering
         \includegraphics[width=.75\textwidth]{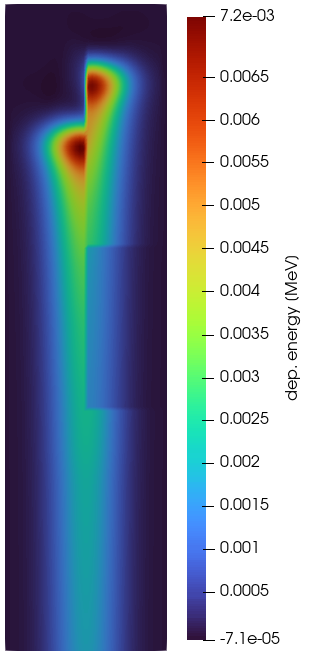}
     \end{subfigure}
          \begin{subfigure}[b]{0.2\textwidth}
         \centering
         \includegraphics[width=.75\textwidth]{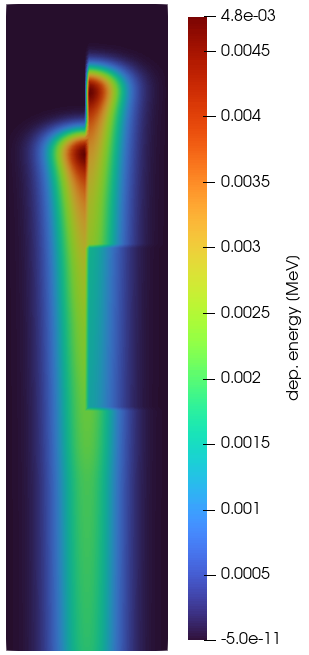}
     \end{subfigure}
          \begin{subfigure}[b]{0.2\textwidth}
         \centering
         \includegraphics[width=.75\textwidth]{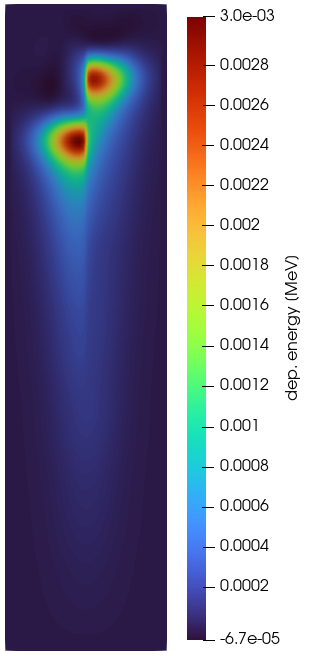}
     \end{subfigure}\\
     \begin{subfigure}[b]{0.2\textwidth}
         \centering
         \includegraphics[width=.75\textwidth]{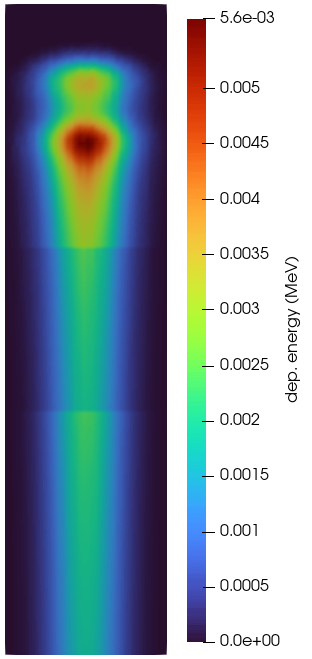}
         \caption{Monte Carlo}
     \end{subfigure}
          \begin{subfigure}[b]{0.2\textwidth}
         \centering
         \includegraphics[width=.75\textwidth]{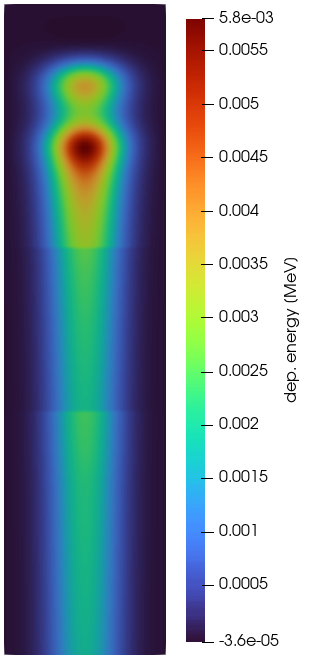}
         \caption{DLRA}
     \end{subfigure}
          \begin{subfigure}[b]{0.2\textwidth}
         \centering
         \includegraphics[width=.75\textwidth]{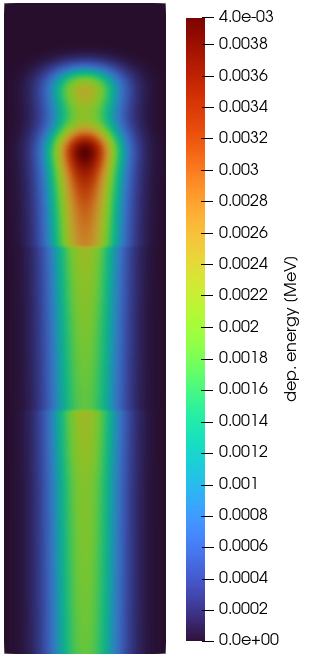}
         \caption{DLRA uncollided}
     \end{subfigure}
          \begin{subfigure}[b]{0.2\textwidth}
         \centering
         \includegraphics[width=.75\textwidth]{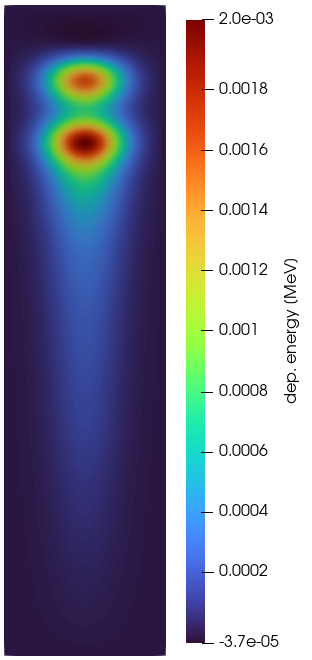}
           \caption{DLRA collided}
     \end{subfigure}
        \caption{Deposited energy of MC vs. DLRA with P$_{95}$, $\vartheta=0.1$ at x=1cm (top) and y=1cm (bottom) in heterogeneous case.}
        \label{fig:HetSurfPlots2}
\end{figure}

Next, we consider a heterogeneous test case which is identical to the previous one apart from a low-density insert with $\mathrm{HU}=-400$ placed between $z=3-5$cm and $y=0-1$cm. 
Again, figures \ref{fig:HetCuts} and \ref{fig:HetSurfPlots2} show that the low-rank solver based on P$_{95}$ is able to approximate the energy deposition well along the central beam axis and at smaller depths. Around the peak, we again observe slightly larger differences to the Monte Carlo reference, however interestingly we do not see a negative energy deposition in this case. This again aligns well with the results from \cite{lathouwers_deterministic_2023}.
In the top row of figure~\ref{fig:HetSurfPlots2} we can however see that in contrast to \cite{lathouwers_deterministic_2023} the diffusive scattering around the peak is not fully captured by the low-rank method. Thus, this more complex set up might also require a finer or more sophisticated angular and/or spatial resolution to further improve accuracy. The higher complexity is also reflected in the chosen ranks, which are higher than in the homogeneous case, especially for the smaller truncation tolerance. The overall agreement to Monte Carlo is slightly lower as well and we see again that the angular resolution has a higher impact than the rank/truncation tolerance (see table \ref{tab:L2diff}).\footnotetext{For NxN matrix with singular values $\sigma_j, j=1,...,N$ the fraction of captured information is defined as $\frac{\sum_{j=1}^{k}\sigma_j^2}{\sum_{j=1}^{N}\sigma_j^2}\cdot100$}
For both test cases, we can see in figure \ref{fig:Ranks} that the ranks at the high initial energy start relatively low. This is likely due to the low scattering cross sections and therefore low contributions of the collided part as well as the shared dominant direction of flight. The ranks then increase as particles start to scatter more and finally decrease again in the end, when the dynamics become almost diffusive for very low energies. The chosen ranks of between 4--5 and 10--22 for the different tolerances are consistently much lower than the degrees of freedom in space (256000 cells) and angle (1444/9025 for P$_{37}$/P$_{95}$), while still managing to capture close to 100\% of the information (see figure \ref{fig:Ranks} (c)), which confirms the validity of a low-rank method.

\section{Conclusions}
In this work, we demonstrated that the computational costs and memory associated with deterministic solvers for a proton transport problem can be significantly reduced using the dynamical low-rank method. The chosen ranks are low and independent of the number of angular modes, making a much higher angular resolution computationally feasible at little additional costs. We show very good agreement with a full-rank solver based on the same physical models and discretizations. Further, our results agree reasonably well with a TOPAS MC reference. Future work could explore the benefits of a higher-order time/energy discretization and extend the framework to not water-equiavalent materials as well as several beam directions.

\printnomenclature

\section*{ACKNOWLEDGEMENTS}
Pia Stammer received funding from the German National Academy of Sciences Leopoldina for the project underlying this article, under grant number LPDS 2024-03.

\bibliographystyle{abbrv}
\bibliography{main}

\appendix

\end{document}